\documentclass[12pt]{article}

\usepackage{a4wide}
\usepackage{amssymb}
\usepackage{amsfonts}
\usepackage{amsmath}
\input xy
\xyoption{arrow} \xyoption{matrix}

\date{}

\newtheorem{proposition}{Proposition}
\newtheorem{theorem}[proposition]{Theorem}

\def\der{\partial }

\def\nFM0{{\nu }_{F,M_0}}
\def\nFN0{{\nu }_{F,N_0}}
\def\nGN0{{\nu }_{G,N_0}}

\def\N0{ {\bf N}_0 }

\def\t{\otimes}

\def\ra{\rightarrow}

\def\Xpm{X^{\pm }}

\def\s{\sigma}

\def\l1{{\lambda}_1}

\def\a{\alpha}
\def\a0{ {\alpha }_0}
\def\a1{ {\alpha }_1}

\def\l{\lambda}

\def\nFGM0{{\nu }_{F,G,M_0}}

\def\nFN0{{\nu}_{F,N_0}}

\def\sm{{\sigma}^m}

\def\sm1{{\sigma}^{-1}}

\def\smtp1{{\sigma}^{-t+1}}

\def\S1{S^{-1}}

\def\Xpm1{X^{\pm 1}_1}

\def\sPM1{{\sigma }^{\pm 1}}
\def\sMP1{{\sigma }^{\mp 1 }}

\def\d{\delta}

\def\di{{\rm d.ind}}

\def\Ytm1{Y^{t-1}}
\def\Yim1{Y^{i-1}}

\def\Aut{{\rm Aut}}

\def\ad{{\rm ad }}

\def\SL2Z{ {\rm SL}_2({\bf Z}) }

\def\Gp1{ G^{1 , 1 } }
\def\P11{ P^{-1 , 1 } }
\def\Pp1{ P^{1 , 1 } }

\def\nCLsr{{}^\nu\kern-2pt {\cal L}^{\sigma , \rho  }}
\def\nP{{}^\nu \kern-2pt P}
\def\nL{{}^\nu\kern-2pt L}
\def\nLL{{}^\nu\kern-2pt \Lambda}
\def\nPsr{{}^\nu\kern-2pt P^{\sigma , \rho  }}
\def\nLsr{{}^\nu\kern-2pt L^{\sigma , \rho  }}
\def\nuCL{{}^\nu\kern-2pt  {\cal L}}
\def\nCLsr{{}^\nu\kern-2pt {\cal L}^{\sigma , \rho  }}
\def\nCL1m{{}^\nu\kern-2pt {\cal L}^{-1 , 1  }}
\def\x1nu{x^\frac{1}{\nu}}
\def\xm1nu{x^{-\frac{1}{\nu}}}

\def\ra{\rightarrow }

\def\nAM0{{\nu }_{{\cal A},M_0}}
\def\nAN0{{\nu }_{{\cal A},N_0}}

\def\ad{ {\rm ad }}

\def\di!{\frac{\der^i}{i!}}
\def\dik!{\frac{\der^k_i}{k!}}
\def\tZ{\otimes_\mathbb{Z}}
\def\Zp{\mathbb{Z}_p}

\begin{document}

\author{V. V. \  Bavula 
}

\title{The ${\rm Jacobian \; Conjecture}_{2n}$ implies the ${\rm Dixmier \;Problem}_n$ }

\maketitle
\begin{abstract}
 Using the {\em inversion formula} for  automorphisms of the Weyl
 algebras with polynomial coefficients and the {\em bound} on its degree \cite{Bav-InvForm} a slightly shorter ({\em algebraic})
 proof is given of the result of  A. Belov-Kanel and M. Kontsevich
 \cite{Bel-Kon05JCDP} that $JC_{2n}$ {\em implies} $DP_n$. No
 originality is claimed.

\end{abstract}



The {\em Weyl } algebra $A_n=A_n(\mathbb{Z})$ is a
$\mathbb{Z}$-algebra generated by $2n$ generators $x_1, \ldots ,
x_{2n}$ subject to the defining relations:
$$ [x_{n+i}, x_j]=\d_{ij}, \;\; [x_i, x_j]=[x_{n+i}, x_{n+j}]=0\;\; {\rm
for\;\; all}\;\; i,j=1, \ldots , n,$$ where $\d_{ij}$ is the
Kronecker delta, $[a,b]:=ab-ba=(\ad a)(b)$. For a ring $R$,
$A_n(R):= R\tZ A_n$ is the Weyl algebra over $R$.
\begin{itemize}
\item {\bf The Jacobian Conjecture $(JC_n)$}: given $\s \in {\rm
End}_{\mathbb{C}-alg} (K[x_1, \ldots , x_n])$ such that  ${\rm
det} (\frac{\der \s (x_i)}{\der x_j})\in K^*:=K\backslash \{ 0\}$
then $\s \in \Aut_{\mathbb{C}} (K[x_1, \ldots , x_n])$. \item {\bf
The
 Dixmier Problem $(DP_n)$}, \cite{Dix}: is a $\mathbb{C}$-algebra endomorphism
 of the Weyl algebra $A_n(\mathbb{C})$ an algebra automorphism?
\end{itemize}

\begin{theorem}\label{i8Nov05}
\cite{Bav-InvForm} {\rm (The Inversion Formula)} For each $\s \in
\Aut_K (A_n(K))$
 and $a\in A_n(K)$,
 $$ \s^{-1}(a)=\sum_{\alpha \in \mathbb{N}^{2n}}\phi_\s
 (\frac{(\der')^\alpha}{\alpha!}a)x^\alpha , $$
 where $x^{\alpha} :=(x_1')^{\alpha_1}\cdots
 (x_{2n}')^{\alpha_{2n}}$,  $(\der')^{\alpha} :=(\der_1')^{\alpha_1}\cdots
 (\der_{2n}')^{\alpha_{2n}}$, $\der_i':= \ad (\s (x_{n+i}))$ and $\der_{n+i}':= -\ad (\s (x_i))$ for $i=1, \ldots , n$,
 $\phi_\s := \phi_{2n}\phi_{2n-1}\cdots \phi_1$ where $\phi_i:= \sum_{k\geq 0}(-1)^i\frac{(\s (x_i))^k}{k!}(\der_i')^k$.
\end{theorem}

{\em Remark}. This result was proved when $K$ is a field of
characteristic zero, but by the {\em Lefschetz principle} it  also
holds for any commutative reduced $\mathbb{Q}$-algebra.

\begin{theorem}\label{degfor}
\cite{Bav-InvForm} Given $\s \in \Aut_K (A_n(K[x_{2n+1}, \ldots ,
x_{2n+m}]))$ where $K$ is a commutative reduced
$\mathbb{Q}$-algebra. Then the degree $\deg \, \s^{-1}\leq (\deg
\, \s)^{2n+m-1}$.
\end{theorem}

\begin{theorem}\label{ThBKK}
\cite{Bel-Kon05JCDP} $JC_{2n}$ $\Rightarrow $ $DP_n$.
\end{theorem}

{\it Proof. }  Let $\s \in {\rm
End}_{\mathbb{C}-alg}(A_n(\mathbb{C}))$.

{\it Step 1.} Let $R$ be a finitely generated (over $\mathbb{Z}$)
$\mathbb{Z}$-subalgebra of $\mathbb{C}$ generated by the
coefficients of the elements $x_i':= \s (x_i)$, $i=1, \ldots ,
2n$. Localizing at finitely many primes $q\in \mathbb{Z}$ one can
assume that the ring $R_p:=R/(p)$ is a domain for all primes $p\gg
0$. Then $\s \in {\rm End}_{R-alg}(A_n(R))$, $x_i'\in A_n(R)=R\tZ
A_n$, and the centre $Z(A_n(R))=R$.

{\it Step 2.}  From this moment on $p\in \mathbb{Z}$ is any (all)
sufficiently big prime number and $\Zp := \mathbb{Z}/(p)$.
\begin{eqnarray*}
 A(p)&:=&  A_n(R)/(p)\simeq R_p\t_{\Zp }A_n(\Zp )\simeq R_p \t_{\Zp }M_{p^n}(\Zp [x_1^p , \ldots , x_{2n}^p] )\\
 &\simeq & M_{p^n}(R_p [x_1^p , \ldots , x_{2n}^p] )=M_{p^n}(C_p)
\end{eqnarray*}
where $x_i^p$ stands for $x_i^p+(p)$, and $M_{p^n}(C_p)$ is a
matrix algebra (of size $p^n$) with coefficients from a polynomial
algebra $C_p:=R_p [x_1^p , \ldots , x_{2n}^p]$ over $R_p$.
 The $\s $ induces an $R_p$-algebra endomorphism $\s_p: A(p)\ra
 A(p)$, $a+(p)\mapsto \s (a)+(p)$.

{\it Step 3.} It follows from the inversion formula (Theorem
\ref{i8Nov05}) and  Theorem \ref{degfor} that

$$ \s\in \Aut_R(A_n(R))\;\; \Leftrightarrow \;\; \s_p\in
\Aut_{R_p}(A(p))\;\; {\rm for \;\; all}\;\; p\gg 0.$$

{\it Step 4.} $\s_p (C_p)\subseteq C_p$ (see \cite{Tsuchi'03}).

{\it Step 5.} Since $A(p)\simeq M_{p^n}(C_p)$, $Z(A(p))=C_p$, and
$\s_p (C_p)\subseteq C_p$, it is obvious that

$$ \s_p\in
\Aut_{R_p}(A(p)) \;\; \Leftrightarrow \;\; \s_p|_{C_p}\in
\Aut_{R_p}(C_p).$$

{\it Step 6. Claim}: $\s_p (C_p)\subseteq C_p$ and $JC_{2n}$ imply
$\s_p|_{C_p}\in \Aut_{R_p}(C_p)$.

{\it Proof of the Claim}. $(i)$.  $(C_p, \{ \cdot , \cdot \} )$ is
a {\em Poisson algebra} where
$$ \{ a+(p), b+(p)\}:= \frac{[a,b]}{p}\;\;  ({\rm mod}\; p)$$
is  the {\em canonical Poisson bracket} on a polynomial algebra in
$2n$ variables (a direct computation, see Lemma 4,
\cite{Bel-Kon05JCDP}) which is obviously $\s_p$-{\em invariant}.

$(ii)$.
\begin{eqnarray*}
 \{\pm 1\} &\ni &  \s_p({\rm det} (\{ x_i^p, x_j^p\})_{1\leq i,j\leq n})={\rm det} (\s_p (\{ x_i^p, x_j^p\}))\\
 &=& {\rm det} ( \{ \s_p(x_i^p), \s_p(x_j^p)\})={\rm det} ( J^T(\{ x_i^p,
 x_j^p\})J)\\
 &=& {\rm det} ( J)^2{\rm det} (\{ x_i^p,
 x_j^p\}))={\rm det} ( J)^2\cdot (\pm 1).
\end{eqnarray*}
where $J:=  (\frac{\der \s (x_i^p)}{\der (x_j^p)})_{1\leq i,j\leq
n}$. Hence, ${\rm det} ( J)\in \{ \pm 1\}$. Only {\em now} we use
the assumption that $JC_{2n}$ holds: which implies $\s_p|_{C_p}\in
\Aut_{R_p}(C_p)$. $\Box$

Department of Pure Mathematics

University of Sheffield

Hicks Building

Sheffield S3 7RH

UK

email: v.bavula@sheffield.ac.uk


\begin{thebibliography}{99}

\bibitem{Bav-InvForm} V. V. Bavula, The inversion formula for automorphisms of the Weyl algebras and  polynomial
algebras, arXiv:math. RA/0512215.

\bibitem{Bel-Kon05JCDP} A. Belov-Kanel and M. Kontsevich, The
Jacobian conjecture is stably equivalent to the Dixmier
Conjecture, arXiv:math. RA/0512171.

\bibitem{Dix} J. Dixmier,  Sur les alg\`{e}bres de Weyl, {\it Bull. Soc. Math.
France}, {\bf 96} (1968), 209--242.

\bibitem{Tsuchi'03} Y. Tsuchimoto, Preliminaries on Dixmier
Conjecture, {\it Mem. Fac. Sci. Kochi Univ. Ser. A Math.}, {\bf
24} (2003), 43--59.
\end{thebibliography}
\end{document}